\documentclass{article}                                     

\input{amssym.def}                                                           
\input{amssym.tex}                                                           
                                                                             
\begin{document}                                                             
\title{Noncommutative geometry as a functor}

\author{Igor  ~Nikolaev
\footnote{Partially supported 
by NSERC.}}


\date{}
 \maketitle


\newtheorem{thm}{Theorem}
\newtheorem{lem}{Lemma}
\newtheorem{dfn}{Definition}
\newtheorem{rmk}{Remark}
\newtheorem{cor}{Corollary}
\newtheorem{prp}{Proposition}
\newtheorem{exm}{Example}
\newtheorem{cnj}{Conjecture}

\newcommand{\N}{{\Bbb N}}
\newcommand{\F}{{\cal F}}
\newcommand{\R}{{\Bbb R}}
\newcommand{\Z}{{\Bbb Z}}
\newcommand{\C}{{\Bbb C}}

\begin{abstract}
In this note the noncommutative geometry is interpreted as 
a functor, whose range is a family of the operator algebras.
Some examples are given and a program is sketched.

\vspace{7mm}

{\it Key words and phrases: functors, operator algebras}

\vspace{5mm}
{\it AMS (MOS) Subj. Class.: 18D;  46L}
\end{abstract}

\section*{Introduction}
A point made in this note is that some  noncommutative spaces  (i.e 
 $C^*$-algebras, Banach or associative algebras)  can  be viewed as a generalized
 homology  in the sense that there exist functors with the range  
in the noncommutative spaces.  The domain of the functors can be any interesting 
category, e.g. the Hausdorff spaces, manifolds,  Riemann surfaces, etc.  
We shall give examples of such  functors.

The above functors have a long history, rather natural and familiar to specialists. 
A foundational example is given
by the Gelfand-Naimark functor, which maps the category of the Hausdorff
spaces to the category of commutative $C^*$-algebras. It was conjectured
by Novikov and proved by Kasparov \cite{Kas1} \&  Mischenko \cite{Mis1}, 
that in many cases the higher signatures of smooth $n$-dimensional manifolds 
are invariants of a certain class of the $C^*$-algebras. The respective 
functor is known as an assembly map. In dynamics, Cuntz \& Krieger \cite{CuKr1}
constructed a functor from the category of topological Markov chains to a category of the 
$C^*$-algebras (Cuntz-Krieger algebras). There are many more examples to add to the list.

As long as a functor is constructed, one can calculate the noncommutative
invariants attached to it. On the face of it, the $C^*$-algebras are a way
 more complex   than the abelian groups.  However,  many important  families of the operator  
algebras have been lately classified in terms of the algebraic $K$-theory 
\cite{ElTo1} and more developments will appear in the future.  The invariants of the $C^*$-algebras 
produce  new (and old) invariants of the objects in the initial  category.
Thus, the problem of an interpretation of the noncommutative invariants   
in terms of the initial category  arises.

In relation to the traditional parts of noncommutative
geometry (e.g. the index theory, cyclic cohomology, quantum groups, etc),
the functorial approach means a switch from a `romantic' to  `pragmatic' 
relationship, in the  sense that the noncommutative spaces become a toolkit 
in the study of the classical spaces.  The problem has two parts:
(i) to map a given  category into a family of the 
noncommutative spaces and  (ii) to prove that the mapping is a functor.  
Note that (ii) is the hardest part of the problem.

The note is organized as follows. In section 1 some examples of
functors with the range in a category of the operator algebras are
considered and  their  noncommutative invariants are analyzed. 
In section 2  draft of a program is sketched.

\bigskip\noindent
{\bf Acknowledgments.} I am grateful to Wolfgang ~Krieger 
for useful discussions and   Ryan ~M. ~Rohm for  
helpful remarks on the first draft of the note.

\section{Three examples}
In this section, some examples of functors with the range in a family
of the noncommutative spaces are given. In the two of three cases, the functors are non-injective.
The list is by no means complete and the reader is encouraged to add examples
of his own.

\subsection{Gelfand and Naimark  functor}
{\bf A.}
This is a foundational example. 
Let $X$ be a locally compact Hausdorff space. By $C(X)$ one understands
a commutative $C^*$-algebra of all functions $f: X\to {\Bbb C}$,
which vanish at infinity. The norm on $C(X)$ is the supremum norm.  
Recall that every point $x\in X$ can be thought of as a linear multiplicative functional
$\hat x: C(X)\to {\Bbb C}$. The Gelfand transform $F: X\to C(X)$
is defined by the formula $x\mapsto f$, where $f\in C(X)$ is  such that
$\hat x(f)=f(x)$.

\medskip\noindent
{\bf B.}
Let $h: X\to Y$ be a continuous map between  the Hausdorff spaces  $X$  and $Y$.
It can be easily shown  that the map
$h_*=F^{-1}\circ h\circ F$ is a homomorphism from the $C^*$-algebra
$C(Y)$ to $C(X)$. In other words, $F$ is a contravariant 
functor from the locally compact Hausdorff spaces to the commutative 
$C^*$-algebras:

\begin{picture}(300,110)(-80,-5)
\put(20,70){\vector(0,-1){35}}
\put(130,70){\vector(0,-1){35}}
\put(98,23){\vector(-1,0){53}}
\put(45,83){\vector(1,0){53}}
\put(5,20){$C(X)$}
\put(118,20){$C(Y)$}
\put(17,80){$X$}
\put(125,80){$Y$}
\put(5,50){$F$}
\put(140,50){$F$}
\put(40,30){\sf homomorphism}
\put(40,90){\sf continuous map}
\end{picture}

\medskip\noindent
{\bf C.} 
Note that $F$ is an injective functor.  The functor $F$ does not
produce new invariants of the Hausdorff spaces,  because of the 
following isomorphism: $K^{alg}(C(X))\cong K^{top}(X)$, where 
$K^{alg}$ and $K^{top}$ are the algebraic and the topological $K$-theory,
respectively.

\subsection{Anosov automorphisms of a two-dimensional torus}
{\bf A.} Let us consider  a non-trivial application of the operator algebras to a problem 
in topology. Recall that an automorphism $\phi: T^2\to T^2$ of the two-dimensional
torus is called {\it Anosov},  if it is given by a matrix 
$A_{\phi}=\left(\small\matrix{a_{11} & a_{12}\cr a_{21} & a_{22}}\right)\in GL(2, {\Bbb Z})$, 
such that $|a_{11}+a_{22}|>2$.  
We wish to construct a functor (an assembly map) 
$\mu: \phi\mapsto {\Bbb A}_{\phi}$, such that for every $h\in Aut~(T^2)$  
the following diagram commutes:

\begin{picture}(300,110)(-80,-5)
\put(20,70){\vector(0,-1){35}}
\put(130,70){\vector(0,-1){35}}
\put(45,23){\vector(1,0){53}}
\put(45,83){\vector(1,0){53}}
\put(5,20){${\Bbb A}_{\phi}\otimes {\cal K}$}
\put(5,55){$\mu$}
\put(140,55){$\mu$}
\put(118,20){${\Bbb A}_{\phi'}\otimes {\cal K},$}
\put(17,80){$\phi$}
\put(117,80){$\phi'=h\circ\phi\circ h^{-1}$}
\put(50,30){\sf isomorphism}
\put(54,90){\sf conjugacy}
\end{picture}

\noindent
where ${\Bbb A}_{\phi}$ is an $AF$-algebra and ${\cal K}$ is the $C^*$-algebra
of compact operators on a Hilbert space.  In other words, if $\phi,\phi'$ are conjugate automorphisms,  
then the $AF$-algebras ${\Bbb A}_{\phi}, {\Bbb A}_{\phi'}$ are stably isomorphic.

\medskip\noindent
{\bf B.}
The map $\mu:\phi\mapsto {\Bbb A}_{\phi}$ is  as follows.
For simplicity, let $a_{11}+a_{22}>2$. 
Note that without loss of generality, one can assume that $a_{ij}\ge 0$ for
a proper basis in the homology group $H_1(T^2; {\Bbb Z})$.  
Consider an $AF$-algebra,  ${\Bbb A}_{\phi}$, given by the following periodic Bratteli diagram:

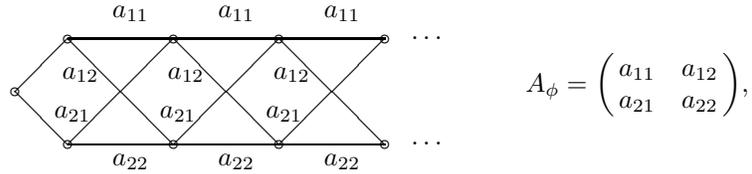
\begin{figure}[here]
\begin{picture}(350,100)(60,0)
\put(100,50){\circle{3}}

\put(120,30){\circle{3}}
\put(120,70){\circle{3}}

\put(160,30){\circle{3}}
\put(160,70){\circle{3}}

\put(200,30){\circle{3}}
\put(200,70){\circle{3}}

\put(240,30){\circle{3}}
\put(240,70){\circle{3}}


\put(100,50){\line(1,1){20}}
\put(100,50){\line(1,-1){20}}

\put(120,30){\line(1,0){40}}
\put(120,30){\line(1,1){40}}
\put(120,70){\line(1,0){40}}
\put(120,70){\line(1,-1){40}}

\put(160,30){\line(1,0){40}}
\put(160,30){\line(1,1){40}}
\put(160,70){\line(1,0){40}}
\put(160,70){\line(1,-1){40}}

\put(200,30){\line(1,0){40}}
\put(200,30){\line(1,1){40}}
\put(200,70){\line(1,0){40}}
\put(200,70){\line(1,-1){40}}

\put(250,30){$\dots$}
\put(250,70){$\dots$}

\put(137,78){$a_{11}$}
\put(177,78){$a_{11}$}
\put(217,78){$a_{11}$}


\put(118,55){$a_{12}$}
\put(158,55){$a_{12}$}
\put(198,55){$a_{12}$}


\put(115,40){$a_{21}$}
\put(155,40){$a_{21}$}
\put(195,40){$a_{21}$}


\put(137,22){$a_{22}$}
\put(177,22){$a_{22}$}
\put(217,22){$a_{22}$}

\put(290,50){
$A_{\phi}=\left(\matrix{a_{11} & a_{12}\cr a_{21} & a_{22}}\right)$,}

\end{picture}

\caption{The $AF$-algebra  ${\Bbb A}_{\phi}$.}
\end{figure}

\noindent
where $a_{ij}$ indicate the multiplicity of the respective edges of the graph. 
We encourage  the reader to verify that $\mu: \phi\mapsto {\Bbb A}_{\phi}$
is a correctly defined function on the set of Anosov automorphisms given by
the hyperbolic matrices with the non-negative entries.
Note that $\mu$ is {\it not} injective, since $\phi$ and all its powers map to
the same $AF$-algebra.

\medskip\noindent
{\bf C.}
Let us show that if $\phi,\phi'\in Aut~(T^2)$ are the conjugate Anosov automorphisms,
then ${\Bbb A}_{\phi},{\Bbb A}_{\phi'}$ are the stably isomorphic
$AF$-algebras.  Indeed, let $\phi'=h\circ\phi\circ h^{-1}$ for
an $h\in Aut~(X)$. Then $A_{\phi'}=TA_{\phi}T^{-1}$ for
a matrix $T\in GL(2, {\Bbb Z})$. Note that $(A_{\phi}')^n=(TA_{\phi}T^{-1})^n=
TA_{\phi}^nT^{-1}$, where $n\in {\Bbb N}$. We shall use the following
criterion (\cite{E}, Theorem 2.3): the $AF$-algebras ${\Bbb A},{\Bbb A}'$
are stably isomorphic if and only if their Bratteli diagrams contain a 
common block of an arbitrary length. Consider the following sequences 
of matrices: $A_{\phi}A_{\phi}\dots A_{\phi}$ and 
$TA_{\phi}A_{\phi}\dots A_{\phi} T^{-1}$,
which mimic the Bratteli diagrams of ${\Bbb A}_{\phi}$ and ${\Bbb A}_{\phi'}$.
Letting  the number of blocks $A_{\phi}$ tend to infinity, we conclude that  
${\Bbb A}_{\phi}\otimes {\cal K}\cong {\Bbb A}_{\phi'}\otimes {\cal K}$.

\medskip\noindent
{\bf D.}
The conjugacy problem for the Anosov automorphisms can now be recast in terms of the $AF$-algebras: 
find invariants of the stable isomorphism classes  of the stationary
$AF$-algebras.  One such invariant is due to Handelman \cite{Han1}. 
Consider an eigenvalue problem for the hyperbolic matrix $A_{\phi}\in GL(2, {\Bbb Z})$:
$A_{\phi}v_A=\lambda_Av_A$,  where $\lambda_A>1$ is the Perron-Frobenius eigenvalue  and 
$v_A=(v_A^{(1)},v_A^{(2)})$ the corresponding eigenvector with the positive entries 
normalized so that $v_A^{(i)}\in K={\Bbb Q}(\lambda_A)$. 
Denote by ${\goth m}={\Bbb Z}v_A^{(1)}+{\Bbb Z}v_A^{(2)}$
 a ${\Bbb Z}$-module in the number field $K$. Recall that the coefficient
ring, $\Lambda$, of module ${\goth m}$ consists of the elements $\alpha\in K$
such that $\alpha {\goth m}\subseteq {\goth m}$. It is known that 
$\Lambda$ is an order in $K$ (i.e. a subring of $K$
containing $1$) and, with no restriction, one can assume that 
${\goth m}\subseteq\Lambda$. It follows  from the definition, that ${\goth m}$
coincides with an ideal, $I$, whose equivalence class in $\Lambda$ we shall denote
by $[I]$. It has been proved by Handelman, that the triple $(\Lambda, [I], K)$ is an arithmetic invariant of the 
stable isomorphism class of ${\Bbb A}_{\phi}$: the ${\Bbb A}_{\phi},{\Bbb A}_{\phi'}$
are stably isomorphic $AF$-algebras if and only if $\Lambda=\Lambda', [I]=[I']$ and $K=K'$. 
It is interesting to compare the operator algebra invariants with those obtained in
\cite{Wal1}.

\medskip\noindent
{\bf E.}
Let $M_{\phi}$ be a mapping torus of the Anosov automorphism $\phi$, i.e.
a three-dimensional manifold $\{T^2\times [0,1] ~|~ (x,0)\mapsto (\phi(x),1) ~\forall x\in T^2\}$.
The $M_{\phi}$ is known as a {\it solvmanifold}, since it is the quotient space of a solvable
Lie group.   It is an easy exercise to show that the homotopy classes of $M_{\phi}$
are bijective with the conjugacy classes of $\phi$. Thus, the noncommutative 
invariant  $(\Lambda, [I], K)$ is a homotopy invariant of $M_{\phi}$.


\subsection{Complex tori and Effros-Shen algebras}
{\bf A.} 
Let us consider an application of the operator algebras to a problem
in conformal geometry. Let $\tau\in {\Bbb H}:=\{z\in {\Bbb C}~|~Im~(z)>0\}$ be 
a complex number.  Recall that the quotient space $E_{\tau}={\Bbb C}/({\Bbb Z}+{\Bbb Z}\tau)$
is called a {\it complex torus}.  It is well-known that 
the complex tori $E_{\tau}, E_{\tau'}$ are isomorphic, whenever 
$\tau'\equiv \tau ~mod~SL(2, {\Bbb Z})$, i.e. $\tau'= {a +b\tau\over c+d\tau}$,  
where $a,b,c,d\in {\Bbb Z}$ and $ad-bc=1$.

\medskip\noindent
{\bf B}. Let $0<\theta< 1$ be an irrational number given by the regular
continued fraction:
$$
\theta=
a_0+{1\over\displaystyle a_1+
{\strut 1\over\displaystyle a_2\displaystyle +\dots}}
$$
By the {\it Effros-Shen algebra} \cite{EfSh1},  
${\Bbb A}_{\theta}$,  one understands  an  $AF$-algebra  given  by the Bratteli diagram:

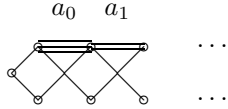
\begin{figure}[here]
\begin{picture}(300,60)(0,0)
\put(110,30){\circle{3}}
\put(120,20){\circle{3}}
\put(140,20){\circle{3}}
\put(160,20){\circle{3}}
\put(120,40){\circle{3}}
\put(140,40){\circle{3}}
\put(160,40){\circle{3}}

\put(110,30){\line(1,1){10}}
\put(110,30){\line(1,-1){10}}
\put(120,42){\line(1,0){20}}
\put(120,40){\line(1,0){20}}
\put(120,38){\line(1,0){20}}
\put(120,40){\line(1,-1){20}}
\put(120,20){\line(1,1){20}}
\put(140,41){\line(1,0){20}}
\put(140,39){\line(1,0){20}}
\put(140,40){\line(1,-1){20}}
\put(140,20){\line(1,1){20}}

\put(180,20){$\dots$}
\put(180,40){$\dots$}

\put(125,52){$a_0$}
\put(145,52){$a_1$}

\end{picture}

\caption{The Effros-Shen algebra  ${\Bbb A}_{\theta}$.}
\end{figure}

\noindent
where $a_i$ indicate the number of edges in the upper row of the diagram.
Recall that ${\Bbb A}_{\theta}, {\Bbb A}_{\theta'}$ are said
to be stably isomorphic if  ${\Bbb A}_{\theta}\otimes {\cal K}
\cong {\Bbb A}_{\theta'}\otimes {\cal  K}$.
It is known that ${\Bbb A}_{\theta}, {\Bbb A}_{\theta'}$ are stably
isomorphic  if $\theta'\equiv \theta~mod~SL(2, {\Bbb Z})$. 
  Comparing the categories of complex tori and Effros-Shen algebras, one cannot fail
to observe   that for the generic objects,  the  corresponding morphism are isomorphic as groups. 
Let us show that the observation has a  ground -- there exists a functor, 
$F$, which makes the following diagram commute:

\begin{picture}(300,110)(-80,-5)
\put(20,70){\vector(0,-1){35}}
\put(130,70){\vector(0,-1){35}}
\put(45,23){\vector(1,0){53}}
\put(45,83){\vector(1,0){53}}
\put(15,20){${\Bbb A}_{\theta}$}
\put(5,55){$F$}
\put(140,55){$F$}
\put(123,20){${\Bbb A}_{\theta'}$}
\put(17,80){$E_{\tau}$}
\put(122,80){$E_{\tau'}$}
\put(40,30){\sf stably isomorphic}
\put(54,90){\sf isomorphic}
\end{picture}

\medskip\noindent
{\bf C.}
To construct the  map $F: E_{\tau}\mapsto {\Bbb A}_{\theta}$,
we shall use a Hubbard-Masur homeomorphism $h: {\Bbb H}\to \Phi_{T^2}$,
where $\Phi_{T^2}$ is the space of measured foliations on the two-torus
\cite{HuMa1}.  Each measured foliation ${\cal F}^{\mu}_{\theta}\in \Phi_{T^2}$ 
looks like a family of the parallel lines of a slope $\theta$ endowed with
an invariant transverse measure  $\mu$ (Fig.3). 
 If $\phi$ is a closed 1-form on $T^2$, then the trajectories of $\phi$
define a measured foliation ${\cal F}^{\mu}_{\theta}\in\Phi_{T^2}$
and vice versa. It is not hard to see that $\mu=\int_{\gamma_1}\phi$
and $\theta=\int_{\gamma_2}\phi /\int_{\gamma_1}\phi$, where 
$\{\gamma_1, \gamma_2\}$ is a basis in $H_1(T^2; {\Bbb Z})$.  
Denote by $\omega_N$ an invariant (N\'eron) differential of the complex torus 
${\Bbb C}/(\omega_1{\Bbb Z}+\omega_2{\Bbb Z})$. It is well known that $\omega_1=\int_{\gamma_1}\omega_N$ and 
$\omega_2=\int_{\gamma_2}\omega_N$.
Let $\pi$ be a projection acting by the formula $(\theta,\mu)\mapsto \theta$. 
The assembly map   $F$  is given by the composition
$F=\pi\circ h$, where $h$ is the Hubbard-Masur homeomorphism. 
In other words, the assembly map $E_{\tau}\mapsto {\Bbb A}_{\theta}$ 
can be written explicitly as: 
$$
E_{\tau}=E_{(\int_{\gamma_2}\omega_N) / (\int_{\gamma_1}\omega_N)}
\buildrel\rm h\over
\longmapsto
{\cal F}^{\int_{\gamma_1}\phi}_{(\int_{\gamma_2}\phi)/(\int_{\gamma_1}\phi)}
\buildrel\rm\pi \over
\longmapsto
{\Bbb A}_{(\int_{\gamma_2}\phi)/(\int_{\gamma_1}\phi)}= {\Bbb A}_{\theta}.
$$

\begin{figure}[here]
\begin{picture}(300,60)(-30,0)

\put(130,10){\line(1,0){40}}
\put(130,10){\line(0,1){40}}
\put(130,50){\line(1,0){40}}
\put(170,10){\line(0,1){40}}

\put(130,40){\line(2,1){20}}
\put(130,30){\line(2,1){40}}
\put(130,20){\line(2,1){40}}
\put(130,10){\line(2,1){40}}

\put(150,10){\line(2,1){20}}

\end{picture}

\caption{The measured foliation ${\cal F}^{\mu}_{\theta}$ on $T^2={\Bbb R}^2/{\Bbb Z}^2$.}
\end{figure}

\medskip\noindent
{\bf D.}
Let us show that the  map $F$ is a covariant functor. 
Indeed, an isomorphism $E_{\tau}\to E_{\tau'}$ is induced by 
an automorphism $\varphi\in Aut~(T^2)$
of the two-torus. Let   
$A_{\varphi}=\left(\small\matrix{a_{11} & a_{12}\cr a_{21} & a_{22}}\right)\in GL(2; {\Bbb Z})$
be a matrix realizing such an automorphism.  From the formulas for $F$,
one gets 
$\tau'=(\int_{c\gamma_1+d\gamma_2}\omega_N)/(\int_{a\gamma_1+b\gamma_2}\omega_N)={c+d\tau\over a+b\tau}$
and 
$\theta'=(\int_{c\gamma_1+d\gamma_2}\phi)/(\int_{a\gamma_1+b\gamma_2}\phi)={c+d\theta\over a+b\theta}$.
Thus, $F$ sends isomorphic complex tori to the stably isomorphic Effros-Shen
algebras. Moreover, the formulas imply that $F$ is a covariant functor. 
Note, that since $F$ contains a  projective map $\pi$, $F$ is not an injective functor.

\medskip\noindent
{\bf E.} Finally, let us consider a noncommutative invariant coming from the
functor $F$. The $E_{CM}$ is said to have a {\it complex multiplication}, 
if the endomorphism ring of the lattice ${\Bbb Z}+{\Bbb Z}\tau$ exceeds ${\Bbb Z}$. 
It is an easy exercise to show  (in view of the explicit fromulas for $F$) that 
$F(E_{CM})={\Bbb A}_{\theta}$, where $\theta$ is a quadratic irrationality. 
In this case the continued fraction of $\theta$  is eventually periodic and 
we let $r$ be the length of the minimal period of $\theta$. Clearly, the integer
$r$ is an invariant of the stable isomorphism class of the $AF$-algebra 
${\Bbb A}_{\theta}$. To interpret the noncommutative invariant $r$ in terms of $E_{CM}$,
recall that $E_{CM}$ is isomorphic to a projective elliptic curve defined 
over a subfield $K={\Bbb Q}(j(E_{CM}))$ of ${\Bbb C}$, where $j(E_{CM})$ is the
$j$-invariant. It is known that the $K$-rational points of $E_{CM}$ make an abelian
group, whose infinite part has rank $R\ge 0$. We conclude by the following
\begin{cnj}
For every elliptic curve with a complex multiplication $R=r-1$. 
\end{cnj}

\section{Sketch of a program}
One can outline a program by indicating:  (i) an object of study,
(ii) a typical  problem and (iii) a set of exercises. A functorial noncommutative
geometry (FNCG) can be described as follows.

\medskip
{\it Object of study.}
The FNCG studies non-trivial functors from a category of the classical objects, ${\goth G}$,
to a category of the noncommutative spaces (operator, Banach or associative algebras), ${\goth A}$.
The functor  can be non-injective.  The category ${\goth A}$ is (possibly) endowed with a good set 
of invariants.

\medskip
{\it  Typical problem.}
The main problem of FNCG is construction of  new invariants of the objects in ${\goth G}$
from the known noncommutative invariants of ${\goth A}$. A reconstruction
of the classical invariants  from the noncommutative invariants
is regarded as a partial solution of the main problem.

\medskip
{\it  Exercises.}
Let ${\goth A}$ be a category of: 

\medskip
(i)  the $UHF$ algebras;

\smallskip
(ii)  the Cuntz-Krieger algebras ${\cal O}_A$ with $det~(A)=\pm 1$.

\medskip\noindent
Find a category ${\goth G}$ corresponding to ${\goth A}$ and solve the typical problem.
(Hint:  for the Cuntz-Krieger algebras of type (ii), the category ${\goth G}$  consists of the homotopy
classes of the torus bundles $M_A$ over $S^1$  with the monodromy given by the matrix $A$.)



\vskip1cm

\textsc{The Fields Institute for Mathematical Sciences, Toronto, ON, Canada,  
E-mail:} {\sf igor.v.nikolaev@gmail.com}

\end{document}